\documentclass[reqno]{amsart}

\usepackage[backref]{hyperref}
\usepackage{mathtools}
\usepackage{tikz}

\newtheorem{theorem}{Theorem}[section]
\newtheorem{lemma}[theorem]{Lemma}

\newtheorem{corollary}[theorem]{Corollary}

\theoremstyle{definition}
\newtheorem{remark}[theorem]{Remark}

\newtheorem{conjecture}[theorem]{Conjecture}
\newtheorem{question}[theorem]{Question}
\newtheorem{notation}[theorem]{Notation}

\usepackage{amssymb}
\renewcommand{\leq}{\leqslant}
\renewcommand{\geq}{\geqslant}
\newcommand{\eps}{\epsilon}

\def\p{\mathfrak p}
\def\q{\mathfrak q}
\def\AGL{\text{AGL}}
\def\PGaL{\text{P$\Gamma$L}}
\def\PGL{\text{PGL}}
\def\PSL{\text{PSL}}

\newcommand\floor[1]{\left\lfloor{#1}\right\rfloor}
\newcommand\ceil[1]{\left\lceil{#1}\right\rceil}

\newcommand{\gen}[1]{\ensuremath{\langle #1\rangle}}

\usepackage{accents}
\newcommand\thickbar[1]{\accentset{\rule{.4em}{.5pt}}{#1}}

\begin{document}

\title{The invariably generating graph of the alternating and symmetric groups}

\author{Daniele Garzoni}
\address{Daniele Garzoni, Universit\`a degli Studi di Padova, Dipartimento di Matematica ``Tullio Levi-Civita''}
\email{daniele.garzoni@phd.unipd.it}
\maketitle

\begin{abstract}
Given a finite group $G$, the invariably generating graph of $G$ is defined as the undirected graph in which the vertices are the nontrivial conjugacy classes of $G$, and two classes are adjacent if and only if they invariably generate $G$. In this paper we study this object for alternating and symmetric groups. The main result of the paper states that, if we remove the isolated vertices from the graph, the resulting graph is connected and has diameter at most $6$.
\end{abstract}

\section{Introduction}
\label{intro}

Given a finite group $G$ and a subset $\{x_1, \ldots, x_t\}$ of $G$, we say that $\{x_1, \ldots, x_t\}$ \textit{invariably generates} $G$ if $\gen{x_1^{g_1}, \ldots, x_t^{g_t}}=G$ for every $g_1, \ldots, g_t \in G$. This concept was introduced by Dixon with motivations from computational Galois theory: see \cite{Dix92} for details. Note that invariable generation can be thought of as a property of conjugacy classes, rather than individual elements.


\subsection{The invariably generating graph} Given a finite group $G$, we define the \textit{invariably generating graph} $\Lambda(G)$ of $G$ as follows. The vertices are the conjugacy classes of $G$ different from $\{1\}$, and two vertices $x^G$ and $y^G$ are adjacent if and only if $\{x,y\}$ invariably generates $G$. The purpose of this paper is to initiate the study of this object for finite (almost) simple groups, and more precisely, for alternating and symmetric groups. It was proved by Kantor--Lubotzky--Shalev \cite{KLS}, and Guralnick--Malle \cite{GM} independently, that finite simple groups are invariably generated by two elements, so that in this case $\Lambda(G)$ is nonempty. It is also known that $S_n$ is invariably generated by two elements if and only if $n\neq 6$ (see \cite[Proposition 4.10]{Tra}). Set
\[
\mathcal P =\left\{\frac{q^d-1}{q-1}: q \,\, \text{is a prime power and} \,\, d\geq 2\right\}.
\]

The first two results of the paper are the following.

\begin{theorem}
\label{main_theorem_1_parte_1}
Assume $n\geq 5$, and let $G\in \{A_n,S_n\}$. Then, $\Lambda(G)$ does not have isolated vertices if and only if $G=A_n$ and $n$ is a prime satisfying $n \notin \mathcal P \cup \{11,23\}$ and $n\equiv -1$ \emph{mod} $12$.
\end{theorem}

\begin{theorem}
\label{main_theorem_1_parte_2}
\begin{enumerate}
    \item Let $(G_i)$ be a sequence of alternating or symmetric groups such that $|G_i|\rightarrow \infty$. Assume that for every $i$, $G_i$ is not alternating of prime degree. Then the number of isolated vertices of $\Lambda(G_i)$ tends to infinity.
    \item Assume $n$ is a prime not contained in $\mathcal P \cup \{11,23\}$. Then the number of isolated vertices of $\Lambda(A_n)$ is at most $2$.
    \end{enumerate}
\end{theorem}

In Theorem \ref{main_theorem_1_parte_1} we assumed $n\geq 5$. We will keep this assumption throughout the paper; see Lemma \ref{lemma_small_degrees} for the remaining cases.

For the sake of clarity, we mention that $\Lambda(A_{11})$ and $\Lambda(A_{23})$ have $5$ and $6$ isolated vertices, respectively (Lemma \ref{M11_M23}).

The only case not addressed in Theorem \ref{main_theorem_1_parte_2} is $G=A_n$ with $n\in \mathcal P$ prime. In Remark \ref{isolated_vertices_n_prime_depends_upon_arithmetic} we will obtain a partial result, dealing with the case $n=(q^d-1)/(q-1)$ with $d\rightarrow \infty$.

Once Theorems \ref{main_theorem_1_parte_1} and \ref{main_theorem_1_parte_2} are proved, one may ask what happens if the isolated vertices are removed from $\Lambda(G)$. With this purpose, we define a graph $\Xi(G)$ which is obtained from $\Lambda(G)$ by removing the isolated vertices. The next result states that, except in case $G=S_6$, this graph is connected with bounded diameter (we already recalled that $S_6$ is not invariably generated by two elements, hence $\Xi(S_6)$ is the null graph).

\begin{theorem}
\label{main_theorem_2}
Assume $n\geq 5$ and let $G\in \{A_n,S_n\}$, with $G\neq S_6$.  Then, $\Xi(G)$ is connected with diameter at most $6$.
\end{theorem}

In many cases we prove better estimates on the diameter.

\begin{theorem}
\label{main_theorem_improvement_diameter}
Assume $n\geq 5$. If $G=S_n$ and $n$ is odd, or if $G=A_n$ and $n$ is even, then $d(\Xi(G))\leq 4$.

If $G=A_n$ and $n\notin \mathcal P$ is prime, then $d(\Xi(G))\leq 3$. 
\end{theorem}

The proofs of Theorems \ref{main_theorem_2} and \ref{main_theorem_improvement_diameter} rely on some recent results, proved in \cite{Jon} and \cite{GMPS}, which classify the primitive subgroups of $S_n$ containing elements having certain cycle types. These results depend on the Classification of Finite Simple Groups.


The upper bound in Theorem \ref{main_theorem_2} is attained, since $d(\Xi(S_8))=6$ (Lemma \ref{lemma_graph_S8}). However, we conjecture that this can happen only for finitely many groups.

\begin{conjecture}
\label{conjecture_diameter_at_most_4}
Let $G \in \{A_n,S_n\}$. If $n$ is sufficiently large then $d(\Xi(G)) \leq 4$.
\end{conjecture}

In Section \ref{subsection_comments_diameter_4} we obtain a partial result towards a proof of this conjecture. It is interesting to observe that in this partial result we use the Prime Number Theorem, but we do not use the CFSG. 

See also Lemma \ref{lower_bound_diameter}, which establishes that $d(\Xi(S_n))\geq 4$ for all primes $n\geq 7$, so the bound in Conjecture \ref{conjecture_diameter_at_most_4}, if true, is attained infinitely often.

Theorem \ref{main_theorem_2} suggests a natural question for all finite simple groups.

\begin{question} Let $G$ be a finite simple group. Is the graph $\Xi(G)$ connected?
\end{question}

\subsection{Some context} The invariably generating graph is the analogue, for invariable generation, of the so called \textit{generating graph} $\Gamma(G)$ of a finite group $G$. This is defined as follows. The vertices are the nonidentity elements of $G$, and two vertices $x$ and $y$ are adjacent if and only if $\gen{x,y}=G$.

Many properties of generation of a finite simple group by two elements can be stated in terms of the generating graph. Guralnick and Kantor \cite{GK} proved that if $G$ is a finite simple group, then for every $1 \neq x \in G$, there exists $y \in G$ such that $\gen{x,y}=G$. Later, Breuer, Guralnick and Kantor \cite{BGK} showed that if $1 \neq x_1, x_2 \in G$, then there exists $y \in G$ such that $\gen{x_1, y}= \gen{x_2,y}=G$. These properties can be stated respectively as follows: $\Gamma(G)$ has no isolated vertices, and it is connected with diameter at most $2$.

Theorems \ref{main_theorem_1_parte_1} and \ref{main_theorem_1_parte_2} say that, for alternating and symmetric groups, the invariably generating graph has quite different properties: it usually has isolated vertices, and the number of isolated vertices usually grows as the order of the group grows.

On the other hand, Theorem \ref{main_theorem_2} says that, if we remove the isolated vertices, we obtain a graph which in some sense shares similarities with the generating graph $\Gamma(G)$.

\subsection{An alternative definition}
In light of Theorem \ref{main_theorem_1_parte_2}, one could ask how the proportion of isolated vertices of $\Lambda(G)$ behaves as $|G|$ tends to infinity. The elementary approach used in the proof of Theorem \ref{main_theorem_1_parte_2} is not sufficient to address this problem. Here however one comment is in order. We have chosen the vertices of $\Lambda(G)$ to be the nontrivial conjugacy classes of $G$. One could define a graph $\Lambda_e(G)$ in which the vertices are the nontrivial \textit{elements} of $G$, and two vertices are adjacent if and only if they invariably generate $G$. Of course, the property of connectedness, and the value of the diameter, are the same in the two graphs. However, when one counts the edges, or the vertices having certain properties, the situation can radically change. Indeed, in $\Lambda_e(G)$ there is a dependence on the size of the conjugacy classes which does not exist in $\Lambda(G)$.

Probabilistic invariable generation has always been considered in terms of elements (cf. \cite{Dix92}, \cite{LP}, \cite{PPR}, \cite{EFG}). Still, we believe it is worth exploring the problem of counting conjugacy classes --- although in this paper we do not address any question of this kind.

\section{Notation and preliminary results}
\label{notation_preliminary_results}
In this section we fix some notation, and we gather some preliminary lemmas and observations that we will use throughout the paper. 

\subsection{Notation} The vertices of our graphs are conjugacy classes of alternating and symmetric groups. We will identify conjugacy classes of $S_n$ with their cycle type, i.e., we will represent conjugacy classes of $S_n$ as partitions of $n$. We now introduce some terminology about partitions.

Let $H \leq S_n$, and let $\mathfrak p$ be a partition of $n$. We will say that \textit{$\p$ belongs to $H$}, or that \textit{$\p$ is contained in $H$}, if $H$ contains elements with cycle type $\p$. If $C$ is the conjugacy class of $S_n$ corresponding to $\p$, this is equivalent to the condition $H \cap C \neq \emptyset$.

Of course, this condition depends only on the $S_n$-conjugacy class of the subgroup $H$. Therefore, in the above terminology we are allowed, if we wish, to replace $H$ by its conjugacy class, and to say that $\p$ belongs to the $S_n$-conjugacy class of $H$.

Let $\p_1$ and $\p_2$ be two partitions of $n$. If $\p_1$ and $\p_2$ both belong to $H$, we will say that $\p_1$ and $\p_2$ \textit{share} $H$.

When $\p \in S_i \times S_{n-i}$, with $1 \leqslant i < n$, we will say that \textit{i is partial sum in $\p$}. This is indeed equivalent to the condition that the integer $i$ can be written as the sum of some parts of $\p$.

In a partition,  $i^m$ will mean $m$ parts of length $i$. Therefore, $(a_1^{m_1}, a_2^{m_2}, \ldots, a_t^{m_t})$ will mean $m_1$ parts of length $a_1$, $m_2$ parts of length $a_2, \ldots$, and $m_t$ parts of length $a_t$.

Occasionally, given a partition $\p$ and a positive integer $i$, we will write $\p^i$ to denote ``the $i$-th power of $\p$'', namely the partition obtained by replacing each part of length $\ell$ by $d$ parts of length $\ell/d$, where $d=(\ell,i)$. Note that if $x\in S_n$ has cycle type $\p$, then $x^i$ has cycle type $\p^i$.

Finally, we define
\[
\mathcal P =\left\{\frac{q^d-1}{q-1}: q \,\, \text{is a prime power and} \,\, d\geq 2\right\}.
\]


\subsection{Maximal overgroups of certain elements} Most of the arguments will rely heavily on the knowledge of the maximal overgroups of certain elements in $S_n$: specifically, cycles, or elements having few orbits in the natural action on $n$ points. The intransitive maximal subgroups are easily determined. For convenience, we now isolate some elementary observations regarding transitive imprimitive subgroups, while then moving to the more difficult case of primitive subgroups. We use some of the language introduced in the previous subsection. The following two lemmas are consequence of \cite[Theorem 2.5]{AACDHL}.

\begin{lemma}
\label{intransitive_two_cycles}
Let $n$ be a natural number, $m$ be a nontrivial divisor of $n$ and $1 \leqslant i <n$. The partition $(a_1, a_2)$ belongs to $S_m \wr S_{n/m}$ if and only if either $m$ divides $a_1$ or $n/m$ divides $a_1$.
\end{lemma}

\begin{proof}
If $(a_1, a_2)$ belongs to $S_m \wr S_{n/m}$, the induced permutation on the blocks has at most two cycles. If it has two cycles, then $m$ divides $a_1$. If it is an $(n/m)$-cycle, then $n/m$ divides $a_1$. The converse implication is proved in the same way.
\end{proof}

\begin{lemma}
\label{intransitive_three_cycles}
Let $n$ be a natural number, $m$ be a nontrivial divisor of $n$ and $1 \leqslant i <n$. The partition $(a_1, a_2, a_3)$ belongs to $S_m \wr S_{n/m}$ if and only if one of the following conditions is satisfied:
\begin{itemize}
\item[(a)] $m$ divides $a_i$ for every $i$.
\item[(b)] $n/m$ divides $a_i$ for every $i$.
\item[(c)] there exist $1 \leqslant t < n/m$ and $i \neq j \in \{1,2,3\}$ such that $a_k=tb_k$ for $k=i,j$, with $b_i+b_j = m$.
\end{itemize}
\end{lemma}

\begin{proof}
Similar to the previous lemma. In case (a), the induced permutation on the blocks has cycle type $(a_1/m, a_2/m, a_3/m)$. In case (b), it is an $(n/m)$-cycle. In case (c), it has cycle type $(t,a_{\ell}/m)$, where $\ell \in \{1,2,3\}\setminus \{i,j\}$.
\end{proof}

We now move to primitive subgroups. Our main tool is a theorem which classifies the primitive subgroups of $S_n$ containing a cycle, and which relies on the CFSG. This should be seen as a generalization of a classical theorem of Jordan (see e.g. \cite[Theorem 13.9]{Wie}) stating that there are no proper primitive subgroups of $S_n$ different from $A_n$ containing a cycle of prime length fixing at least $3$ points. Since we will apply this result several times, for convenience we report here the statement.

\begin{theorem} \cite{Jon}
\label{jones_theorem}
Let $G$ be a primitive permutation group of finite degree $n$, not containing the alternating group $A_n$. Suppose that $G$ contains a cycle fixing $k$ points, where $0 \leqslant k \leqslant n-2$. Then one of the following holds:
\begin{enumerate}
\item $k = 0$ and either
\begin{itemize}
\item[(a)] $C_p \leqslant G \leqslant \AGL_1(p)$ with $n = p$ prime, or
\item[(b)] $\PGL_d(q) \leqslant G \leqslant P\Gamma L_d(q)$ with $n = (q^d-1)/(q-1)$ and $d \geqslant 2$ for some prime
power $q$, or
\item[(c)] $G = PSL_2(11), M_{11}$ or $M_{23}$ with $n = 11, 11$ or $23$ respectively.
\end{itemize}
\item $k = 1$ and either
\begin{itemize}
\item[(a)] $\AGL_d(q) \leqslant G \leqslant A\Gamma L_d(q)$ with $n = q^d$ and $d \geqslant 1$ for some prime power $q$, or
\item[(b)] $G = \PSL_2(p)$ or $\PGL_2(p)$ with $n = p + 1$ for some prime $p \geqslant 5$, or
\item[(c)] $G = M_{11}$, $M_{12}$ or $M_{24}$ with $n = 12, 12$ or $24$ respectively.
\end{itemize}
\item $k = 2$ and $\PGL_2(q) \leqslant G \leqslant P\Gamma L_2(q)$ with $n = q + 1$ for some prime power $q$.
\end{enumerate}
\end{theorem}

Note that the statement implies that there are no proper primitive subgroups of $S_n$ different from $A_n$ containing a cycle fixing at least $3$ points, generalizing indeed Jordan's theorem. We note the following immediate consequence.

\begin{corollary}
\label{cycle_fixing_3_points_power}
Assume $x \in S_n$ is such that a suitable power of $x$ is a nontrivial cycle fixing at least $3$ points. Then, $x$ does not lie in proper primitive subgroups of $S_n$ different from $A_n$.
\end{corollary}

We also mention that we will make essential use of the main result from \cite{GMPS}, which classifies the primitive subgroups of $S_n$ containing an element having at most $4$ cycles.

We will shortly apply the previous results to certain elements (or partitions) of particular interest to us.

\subsection{Conjugacy classes of $A_n$}
If a partition of $n$ is made of distinct odd parts, then the corresponding $S_n$-conjugacy class splits into two $A_n$-conjugacy classes (and viceversa), giving rise to two vertices of $\Lambda(A_n)$. Often, this does not represent a serious change: the following technical lemmas give conditions under which the two vertices may be essentially thought of as a unique vertex.

\begin{lemma}
\label{maximal_subgroups_conjugacy}
Let $x \in A_n$, and let $H \leqslant S_n$ be such that $H^{A_n}=H^{S_n}$. Then, $H$ contains elements belonging to $x^{A_n}$ if and only if it contains elements belonging to $(x')^{A_n}$ for every $x'\in x^{S_n}$.
\end{lemma}

\begin{proof}
Assume $x^z \in H$ with $z \in A_n$, and assume $x'=x^{zg}$ with $g \in S_n$. Then $x'=x^{zg} \in H^g=H^h$ for some $h \in A_n$ by hypothesis, hence $(x')^{h^{-1}} \in H$. This concludes the proof.
\end{proof}

\begin{notation}
\label{notation_unico_vertice}
Assume $x,y \in A_n$ are such that $(x')^{A_n}$ is adjacent to $(y')^{A_n}$ in $\Lambda(A_n)$ for any $x' \in x^{S_n}$ and for any $y' \in y^{S_n}$. Under these assumptions, we say with slight abuse of notation that $x^{S_n}$ is adjacent to $y^{S_n}$ in $\Lambda(A_n)$.
\end{notation}

This notation will be convenient, as we will represent $S_n$-conjugacy classes as partitions, and we will be allowed to say that ``a partition $\mathfrak q$ is adjacent to a partition $\p$'', rather than ``any $A_n$-conjugacy class of elements with cycle type $\mathfrak q$ is adjacent to any $A_n$-conjugacy classes of elements with cycle type $\p$''. We will now see that in many cases the assumption of Notation \ref{notation_unico_vertice} is satisfied.

\begin{lemma}
\label{unico_vertice}
Let $x\in A_n$, and assume $\gen x^{A_n}=\gen x^{S_n}$. Let $y \in A_n$. Then, $x^{A_n}$ is adjacent in $\Lambda(A_n)$ to $y^{A_n}$ if and only if $x^{S_n}$ and $y^{S_n}$ are adjacent (in the terminology of Notation \ref{notation_unico_vertice}).
\end{lemma}

\begin{proof}
Choose $x' \in x^{S_n}$. We show first that if $y^{A_n}$ is not adjacent to $x^{A_n}$, then $y^{A_n}$ is not adjacent to $(x')^{A_n}$. By assumption we can write $x'=(x^i)^z$ for some integer $i$ and for some $z \in A_n$. Again by assumption, $\{x^{g_1}, y^{g_2}\} \subseteq H$ for some $g_1, g_2 \in A_n$, and for some proper subgroup $H$ of $A_n$. Then $(x')^{z^{-1}g_1} = (x^i)^{g_1} = (x^{g_1})^i \in H$, whence $\{(x')^{z^{-1}g_1}, y^{g_2}\} \subseteq H$ and $y^{A_n}$ is not adjacent to $(x')^{A_n}$, as required. 

Assume now $y'\in y^{S_n}$. We show that if $y^{A_n}$ is not adjacent to $x^{A_n}$, then $(y')^{A_n}$ is not adjacent to $x^{A_n}$. This will conclude the proof. We may assume $y^{A_n} \neq (y')^{A_n}$ and $x^{A_n}\neq x^{S_n}$, otherwise the statement is easy. Choose $x' \in x^{S_n}\setminus x^{A_n}$. By the previous paragraph, if $y^{A_n}$ is not adjacent to $x^{A_n}$ then it is not adjacent to $(x')^{A_n}$, so that $\{x^{g_1}, y^{g_2}\} \subseteq H$ for some proper subgroup $H$ of $A_n$ and for some $g_1 \in S_n\setminus A_n$ and $g_2 \in A_n$. Then $\{x, y^{g_2g_1^{-1}}\} \subseteq H^{g_1^{-1}}$, and since $y^{g_2g_1^{-1}}$ is $A_n$-conjugate of $y'$, the proof is concluded.
\end{proof}

We now apply the previous considerations to certain elements and subgroups of $S_n$.

\begin{lemma}
\label{intransitive_imprimitive}
Let $H$ be a maximal subgroup of $S_n$ which is either intransitive, or transitive and imprimitive. Then $H^{A_n}=H^{S_n}$.
\end{lemma}

\begin{proof}
We may assume $n \geq 3$. We have that $H^{A_n}=H^{S_n}$ if and only if $\text{N}_{S_n}(H) \nleqslant A_n$. In our case, $H \nleqslant A_n$, since every maximal intransitive or imprimitive subgroup of $S_n$ contains transpositions.
\end{proof}

\begin{lemma}
\label{only_intrasitive_and_imprimitive_matter}
Assume $x \in A_n$ belongs to no proper primitive subgroup of $S_n$ different from $A_n$. Let $y \in A_n$. Then, $x^{A_n}$ and $y^{A_n}$ are adjacent in $\Lambda(A_n)$ if and only if $x^{S_n}$ and $y^{S_n}$ are adjacent (in the terminology of Notation \ref{notation_unico_vertice}).
\end{lemma}

\begin{proof} By assumption, $x$ lies in no proper primitive subgroups different from $A_n$. Therefore, $x^{A_n}$ is adjacent to $y^{A_n}$ if and only if, for every $g_1,g_2 \in A_n$, $\gen{x^{g_1}, y^{g_2}}$ is primitive or, equivalently, is not contained in intransitive or imprimitive maximal subgroups. By Lemmas \ref{maximal_subgroups_conjugacy} and \ref{intransitive_imprimitive}, this condition depends only on the cycle type of the elements, rather than on their $A_n$-conjugacy class. The lemma follows.
\end{proof}

\begin{lemma}
\label{i_coprime_n}
Let $2 \leq i \leq n/2$ be coprime with $n$. Then, $(i,n-i)$ does not belong to proper primitive subgroups different from $A_n$, and does not belong to transitive imprimitive subgroups.
\end{lemma}

\begin{proof}
The fact that $(i,n-i)$ does not lie in imprimitive subgroups follows from Lemma \ref{intransitive_two_cycles}. Regarding primitive subgroups, $(i,n-i)^{n-i}=(i, 1^{n-i}) \neq (1^n)$, and $n-i \geq 3$ (the conditions on $i$ imply $n \geq 5$), hence the statement follows from Corollary \ref{cycle_fixing_3_points_power}.
\end{proof}

Theorem \ref{jones_theorem} can be used to generalize the previous lemma to the case $i=1$: one just needs to take care of some specific examples of primitive subgroups. The following easy lemma will be used with this purpose.

\begin{lemma}
\label{affine_and_projective} Let $q$ be a prime power and $d$ be a positive integer.
\begin{enumerate}
    \item An element of $\AGL_d(q)$, in the natural action on $q^d$ points, either is a derangement, or fixes a number of points equal to $q^s$ for some $0 \leq s \leq d$.
    \item Assume $g\in P\Gamma L_2(q)$ fixes at least $3$ points in the natural action on $q+1$ points. Then, $g$ fixes a number of points having the same parity of $q+1$. Moreover, if $g\in \PGL_2(q)$ then $g=1$.
\end{enumerate}
\end{lemma}

\begin{proof}
(1) If $g \in \AGL_d(q)$ fixes some point, we may assume that it fixes $0$. Hence, $g \in \text{GL}_d(q)$. Now just observe that the set of fixed points of an element of $\text{GL}_d(q)$ is an $\mathbf F_q$-subspace of $\mathbf{F}_q^d$.

(2) Consider $\PGaL_2(q)$ acting (on the right) on the set $\Omega$ of $1$-dimensional subspaces of $\mathbf F_q^2$. Write $q=p^r$ with $p$ prime. For $\phi \in \text{Gal}(\mathbf{F}_{q}/\mathbf{F}_{p})$, denote by $f_\phi$ the permutation of $\Omega$ induced by the mapping $(\lambda_1,\lambda_2)\mapsto (\lambda_1^\phi, \lambda_2^\phi)$ of $\mathbf F_q^2$. Then, we may express each element $g\in \PGaL_2(q)$ as $g=xf_\phi$, where $x\in \text{PGL}_2(q)$ and $\phi \in \text{Gal}(\mathbf{F}_{q}/\mathbf{F}_{p})$.

Note that $\PGaL_2(q)$ is $3$-transitive on $\Omega$. Hence, if $g=x f_\phi \in \PGaL_2(q)$ fixes at least $3$ points, we may assume that it fixes $\gen{(1,0)}$, $\gen{(0,1)}$ and $\gen{(1,1)}$. It follows that $x=1$ and $g=f_\phi$. Then $g$ fixes $p^\ell+1$ points, where $\ell$ is a divisor of $r$; in particular, it fixes a number of points having the same parity of $q+1$. Moreover, if $g \in \PGL_2(q)$ then $f_\phi=1$. The lemma is proved.
\end{proof}

\begin{lemma}
\label{n_cycle}
Assume $n$ is either an odd prime, or $n\equiv 3$ mod $4$. Assume $x \in A_n$ is an $n$-cycle. Then, $\gen x^{A_n}=\gen x^{S_n}$.
\end{lemma}

\begin{proof}
If $n\geq 3$ is an odd integer, then an $n$-cycle is normalized by elements having cycle type $(1,2^{(n-1)/2})$. If $n$ is an odd prime, then an $n$-cycle is normalized by an $(n-1)$-cycle. In particular, if $x$ is as in the statement, then $\text N_{S_n}(\gen x) \nleqslant A_n$, from which $\gen x^{A_n}=\gen x^{S_n}$.
\end{proof}

We deduce another consequence of the previous lemmas.

\begin{lemma}
\label{automorphism_doesnt_do_anything}
Assume $n \geq 4$, and let $x \in A_n$. Then, for every $g \in S_n$, $x^{A_n}$ and $(x^g)^{A_n}$ are not adjacent in $\Lambda(A_n)$.
\end{lemma}

\begin{proof}
By Lemma \ref{intransitive_imprimitive}, we may assume that $x$ is not contained in intransitive or imprimitive subgroups. It follows that $n$ is prime and $x$ is an $n$-cycle. The statement follows then by Lemma \ref{n_cycle}.
\end{proof}

This lemma suggests a natural question for all finite simple groups.

\begin{question} Let $G$ be a finite simple group. Let $x \in G$, and let $\sigma \in \text{Aut}(G)$. Is it possible that $\{x,x^{\sigma}\}$ invariably generates $G$?
\end{question}

It is easy to check that the statement is true for $G=A_6$, in which case $S_6$ has index two in $\text{Aut}(G)$. This, together with Lemma \ref{automorphism_doesnt_do_anything}, implies that the question has a negative answer for $G=A_n$.
\subsection{Small degrees} The main theorems are stated for $n\geq 5$. For completeness, we address here the cases of degree $n\leq 4$. We assume $n\geq 3$ in order to avoid trivialities. 

\begin{lemma}
\label{lemma_small_degrees}
Let $3\leq n \leq 4$, and assume $G\in \{A_n,S_n\}$. Then, $\Lambda(G)$ has isolated vertices if and only if $G=S_4$. Moreover, $\Xi(G)$ is connected with diameter at most $2$.
\end{lemma}

\begin{proof}
This is an easy check. The graphs $\Xi(A_3), \Xi(S_3)$ and $\Xi(S_4)$ have diameter $1$, while $\Xi(A_4)$ has diameter $2$: the two classes of $3$-cycles are connected by a path of length $2$ passing through the class $(2^2)$. In $\Lambda(S_4)$, the vertices corresponding to $(1^2,2)$ and $(2^2)$ are isolated.
\end{proof}

\section{Proofs} In this section we prove the theorems stated in the introduction.

\subsection{Proof of Theorems \ref{main_theorem_1_parte_1} and \ref{main_theorem_1_parte_2}}
\label{proof_main_theorem_1}
We begin with a lemma which proves Theorem \ref{main_theorem_1_parte_2}(2) and the ``if'' part of Theorem \ref{main_theorem_1_parte_1}.

\begin{lemma}
\label{isolated_vertices_An_n_prime}
Let $n \geq 5$ be an odd prime. Then, $(1,2^{(n-1)/2})$ and $(1,3^{(n-1)/3})$ are isolated in $\Lambda(A_n)$ (when they make sense and are even permutations). If $n\notin \mathcal P \cup \{11,23\}$, then there are no other isolated vertices in $\Lambda(A_n)$.
\end{lemma}

\begin{proof} The two mentioned vertices are not adjacent to $(n)$ (this terminology makes sense by Lemmas \ref{unico_vertice}, \ref{n_cycle} and Notation \ref{notation_unico_vertice}) because they are contained in $\AGL_1(n)$. Therefore, they might be adjacent only to some $A_n$-conjugacy class with cycle type $(a_1, \ldots, a_t)$, with $t \geqslant 3$. Every $1 \leqslant i \leqslant (n-1)/2$ is a partial sum in $(1,2^{(n-1)/2})$, so this vertex is isolated (recall Lemmas \ref{intransitive_imprimitive} and \ref{maximal_subgroups_conjugacy}), and the values of $i$ that are not partial sum in $(1,3^{(n-1)/3})$ are exactly those that satisfy $i \equiv 2$ mod $3$. If $a_1, a_2 \equiv 2$ mod $3$ then $a_1+a_2 \equiv 1$ mod $3$, so one among $a_1, a_2$ and $a_1+a_2$ is partial sum. We conclude that $(1,3^{(n-1)/3})$ is indeed isolated.

Assume now $n\notin \mathcal P \cup \{11,23\}$. Let $x^{A_n}$ be a vertex of $\Lambda(A_n)$ different from the two above: we want to show it is not isolated. If $x \notin \AGL_1(n)$, we deduce from Theorem \ref{jones_theorem} that $x^{A_n}$ is adjacent to $(n)$. On the other hand, if $x \in \AGL_1(n)$, then either it is an $n$-cycle, or it has cycle type $(1,t^{(n-1)/t})$, with $t \geqslant 4$. As just remarked, a class of $n$-cycles is not isolated. If $t \geqslant 5$, by Corollary \ref{cycle_fixing_3_points_power} $x^{A_n}$ is adjacent to $(2^2,n-4)$. If $t=4$, for the same reason $x^{A_n}$ is adjacent to $(3^2,n-6)$.
\end{proof}

For the sake of clarity, we deal with the cases $n=11,23$.

\begin{lemma}
\label{M11_M23}
The graphs $\Lambda(A_{11})$ and $\Lambda(A_{23})$ have $5$ and $6$ isolated vertices, respectively.
\end{lemma}

\begin{proof}
Recall Theorem \ref{jones_theorem}. Let $n\in \{11,23\}$. Let $1 \neq x\in A_n$. If $x$ does not belong to $M_{11}$ (for $n=11$) and $M_{23}$ (for $n=23$), the same argument as in the previous lemma shows that $x^{A_n}$ is adjacent to $(n)$ (note that $\text{PSL}_2(11)$ can be embedded in $M_{11}$). Inspection (using for instance GAP) shows that if $n=11$ and $1 \neq x\in M_{11}$, then $x^{A_n}$ belongs to
\[
\{(1^3,2^4), (1^3,4^2), (1^2,3^3), (2,3,6), (1,2,8)\}.
\]
These are not adjacent to $(n)$ because of $M_{11}$. By looking at partial sums, we see that these are all isolated. If $n=23$ and $1 \neq x\in M_{23}$, then either $x^{A_n}$ is adjacent to one between $(2^2,19)$ and $(3^2, 17)$, or $x^{A_n}$ belongs to
\[
\{(1^5,3^6), (1^3,5^4), (1^7,2^8), (1,2^2,3^2,6^2), (1^3,2^2,4^4), (1,2,4,8^2)\}.
\]
These are all isolated for the same reason as above.
\end{proof}

We are now ready to prove Theorems \ref{main_theorem_1_parte_1} and \ref{main_theorem_1_parte_2}.

\begin{proof}[Proof of Theorems \ref{main_theorem_1_parte_1} and \ref{main_theorem_1_parte_2}]
In the proof, we will use Lemmas \ref{maximal_subgroups_conjugacy}, \ref{unico_vertice} and \ref{intransitive_imprimitive} with no further mention. Recall also Notation \ref{notation_unico_vertice}.

We begin with Theorem \ref{main_theorem_1_parte_2}. Item (2) follows from Lemma \ref{isolated_vertices_An_n_prime}, hence we focus on item (1). Let $G=A_n$ or $S_n$, with $n$ nonprime if $G=A_n$. For every $n$, we will define $I_n$, a subset of the set of isolated vertices of $\Lambda(G)$, whose cardinality, as $n \rightarrow \infty$, goes to infinity. This will prove Theorem \ref{main_theorem_1_parte_2}. For every $m$, denote by $p_m^E$ the set of all partitions of $m$ different from $(1^m)$ which correspond to even permutation.

We first assume that, whenever $n$ is odd, then $G=S_n$. Define $I_n$ as the set of all partitions of $n$ of the form $(1^{\floor{n/2}}, \p)$, where $\p$ is whatsoever partition of $p_{\ceil{n/2}}^E$.

Let $\q \in I_n$. Every $1 \leqslant i \leqslant n/2$ is partial sum, so $\q$ is not adjacent to classes of elements with at least two cycles. If $n$ is even, then $\q$ is not adjacent to $(n)$ since $n/2$ is partial sum in $\q$, hence $\q \in S_{n/2} \wr S_2$. If $n$ is odd, then $G=S_n$, and $\q$ is not adjacent to $(n)$ since $\q$ corresponds to even permutations. 

Therefore, $I_n$ consist of isolated vertices. Clearly, the size of $I_n$ goes to infinity as $n \rightarrow \infty$.

Assume now that $G=A_n$ and $n$ is odd nonprime. Fix $n$, and let $p_n$ be the smallest prime divisor of $n$. Define $I_n$ as the set of all partitions of $n$ of the form $(1^{n(p_n-1)/p_n}, \p)$, where $\p$ is any partition of $p_{n/{p_n}}^E$.

If $\q \in I_n$ then $\q \in S_{n/{p_n}} \wr S_{p_n}$: the first $p_n-1$ blocks are fixed pointwise. Therefore, $\q$ is not adjacent to $(n)$. Moreover, every $1 \leqslant i \leqslant n/2$ is partial sum, so $\q$ is not adjacent to classes of elements having at least $2$ cycles. It follows that $\q$ is isolated. Note now that $n/p_n \geqslant \sqrt{n}$, hence the size of $I_n$ goes to infinity as $n\rightarrow \infty$. This concludes the proof of Theorem \ref{main_theorem_1_parte_2}.

We now move to Theorem \ref{main_theorem_1_parte_1}. Note that Lemma \ref{isolated_vertices_An_n_prime} proves the ``if'' part. We now prove the ``only if'' part. Assume first $n$ is prime and $G=A_n$. If $n \not\equiv -1$ mod $12$, by Lemma \ref{isolated_vertices_An_n_prime} there are isolated vertices in $\Lambda(A_n)$. The cases $n=11,23$ have been considered in Lemma \ref{M11_M23}. If $n=(q^d-1)/(q-1)$, let $x$ be any involution lying in a subgroup of $S_n$ conjugate to $\PGaL_d(q)$. The fact that $n$ is odd implies that every $1 \leq i \leq n/2$ is partial sum in (the cycle type of) $x$, hence $x^{A_n}$ might only be adjacent to a class of $n$-cycles. However, this does not happen because of the containment in $\PGaL_d(q)$. Therefore $x^{A_n}$ is isolated.

The case $G=A_n$ with $n$ prime is therefore proved. For the remaining cases, we apply what proved for Theorem \ref{main_theorem_1_parte_2}. We define $\q=(1^{n-3},3)$. This partition belongs to $I_n$, as defined in this proof, hence it is an isolated vertex of $\Lambda(G)$. This concludes the proof.
\end{proof}

\begin{remark}
\label{isolated_vertices_n_prime_depends_upon_arithmetic}
The unique case not discussed in Theorem \ref{main_theorem_1_parte_2} is the case $G=A_n$ and $n\in \mathcal P$ prime. We obtain here the following partial result: if $n_i=(q^{d_i}-1)/(q-1)$ is prime and $d_i\rightarrow \infty$, then the number of isolated vertices of $\Lambda(A_{n_i})$ tends to infinity. Note that since $n_i$ is prime, $d_i$ must be prime. To establish whether infinitely many such primes (i.e., primes $d$ such that $(q^d-1)/(q-1)$ is prime for some prime power $q$) do actually exist, however, is a hard open problem in number theory: see for instance \cite{BLS}.

Assume first $n=(q^d-1)/(q-1)$ with $q$ odd (and $d$ odd). The action of $\PGL_d(q)$ on the $1$-dimensional subspaces of $\mathbf F_q^d$ gives an embedding $\PGL_d(q) < S_n$. For every $1 \leq r < d/2$, let $x=x(r)$ be a diagonal matrix of $\text{GL}_d(q)$ with $1$'s and $-1$'s on the diagonal, and assume the number of $-1$'s is $r$. Let $\thickbar x$ denote the image of $x$ in $\PGL_d(q) < S_n$. Then, $\thickbar x$ has $(q^r+q^{d-r}-2)/(q-1)$ fixed points. In particular, any two distinct $1 \leq r < d/2$ give rise to elements of $S_n$ which have a different number of fixed points, and which therefore belong to different $S_n$-conjugacy classes. It is easy to check that the element $\thickbar x$ arising in this way belongs to $A_n$. The number of possibilities for $r$ in order to obtain such an element is $(d-1)/2$. As remarked in previous proof, $(\thickbar x)^{A_n}$ is isolated in $\Lambda(A_n)$; therefore the number of isolated vertices of $\Lambda(A_n)$ is at least $(d-1)/2$.

Assume now $q$ is even. For every $1\leq \ell \leq d/2$, consider a unipotent element $x=x(\ell)$ of $\text{SL}_d(q)$ with $\ell$ Jordan blocks of size $2$, and with the other Jordan blocks of size $1$. This is an involution of $\text{SL}_d(q)$. Denote again by $\thickbar x$ the image of $x$ in $\PSL_d(q)< S_n$. Then $\thickbar x$ has $(q^{d-\ell}-1)/(q-1)$ fixed points, hence any two distinct $1 \leq \ell \leq d/2$ give rise to elements belonging to different $S_n$-conjugacy classes. Moreover $\thickbar x \in A_n$ (unless $(q,d)=(2,2)$, but recall we are assuming $n\geq 5$), and $(\thickbar x)^{A_n}$ is isolated in $\Lambda(A_n)$. Hence the number of isolated vertices of $\Lambda(A_n)$ is at least $(d-1)/2$.
\end{remark}

\subsection{Proof of Theorems \ref{main_theorem_2} and \ref{main_theorem_improvement_diameter}} In this subsection we prove Theorems \ref{main_theorem_2} and \ref{main_theorem_improvement_diameter}. 

One brief comment about the terminology we will adopt. The proofs will begin with a sentence of the type ``Let $\p$ be a vertex of $\Xi(G)$'', without any preliminary consideration showing that $\Xi(G)$ is not the null graph. However, along the proof suitable edges will be exhibited in $\Lambda(G)$, so that the initial choice of $\p$ will be licit. (In other words we are saying that, although we will not state it explicitly, in the proofs it will be shown that the groups are invariably generated by two elements.)

We begin by proving Theorem \ref{main_theorem_improvement_diameter} (in two separate results).

\begin{theorem}
\label{lemma_better_estimate_diameter}
Assume $n\geq 5$ is a prime and $n\notin \mathcal P$. Then $\Xi(A_n)$ is connected and $d(\Xi(A_n))\leq 3$.
\end{theorem}

\begin{proof}
Let $x \in A_n$, and assume $x^{A_n}$ is a vertex of $\Xi(A_n)$. If $n \neq 11, 23$, in the proof of Lemma \ref{isolated_vertices_An_n_prime} we showed that $x^{A_n}$ is adjacent to one among $(n), (2^2,n-4)$ and $(3^2,n-6)$ (recall Lemma \ref{n_cycle}). These vertices are pairwise adjacent by Corollary \ref{cycle_fixing_3_points_power}, hence we have indeed $d(\Xi(A_n)) \leq 3$. If $n=11$ or $23$, we observed in the proof of Lemma \ref{M11_M23} that a class inside $M_{11}$ or $M_{23}$ is either isolated, or adjacent to one between $(2^2,n-4)$ and $(3^2,n-6)$. Therefore, by Corollary \ref{cycle_fixing_3_points_power} also in this case $d(\Xi(A_n)) \leq 3$.
\end{proof}

For a later use, we point out that the same estimate holds for the groups $A_{13}$ and $A_{17}$.

Assume $n=13=(3^3-1)/(3-1)$. Inspection shows that every class lying in $\PGL_3(3)$ is either isolated or adjacent to $(3^3, 7)$, hence $d(\Xi(A_{13})) \leq 3$ by the same argument as in the proof of the previous theorem. Assume finally $n=17=16+1$. Here every class lying in $\PGaL_2(16)$ is either isolated or adjacent to one between $(3^2,11)$ and $(4^2,9)$. It is easy to deduce that $d(\Xi(A_{17})) \leq 3$.

\begin{theorem}
\label{symemtric_odd_degree_and_alternating_even_degree}
Let $n\geq 5$ be an integer.

\begin{enumerate}
\item If $n$ is odd then $\Xi(S_n)$ is connected and $d(\Xi(S_n)) \leqslant 4$.
\item If $n$ is even then $\Xi(A_n)$ is connected and $d(\Xi(A_n)) \leqslant 4$.
\end{enumerate}

\end{theorem}

\begin{proof}
(1) Let $\p$ be a vertex of $\Xi(S_n)$. Assume $\p$ is not adjacent to $(i,n-i)$ for any $2 \leqslant i \leqslant n/2$ coprime with $n$. Then, by Lemma \ref{i_coprime_n} every such $i$ is partial sum in $\p$. Since $\p$ is a vertex of $\Xi(S_n)$, $\p$ will be adjacent to some vertex $\q$. Necessarily, no $2 \leqslant i \leqslant n/2$ coprime with $n$ is partial sum in $\q$. Hence, again by Lemma \ref{i_coprime_n}, $\q$ is adjacent to $(i,n-i)$ for every such $i$.

Now note that the $(i,n-i)$'s, with $i$ as above, are pairwise adjacent. From this it follows that, for any fixed $2 \leq i \leq n/2$ coprime with $n$, any vertex $\p$ of $\Xi(S_n)$ has distance at most $2$ from $(i,n-i)$. This concludes the proof.

(2) The statement for $n=6$ can be checked explicitly, hence we assume $n\geq 8$. Then the proof is identical to (1). Recall Lemmas \ref{maximal_subgroups_conjugacy} to  \ref{i_coprime_n}, and Notation \ref{notation_unico_vertice}.
\end{proof}

Now we move to the proof of the general case, i.e., Theorem \ref{main_theorem_2}. We are left with symmetric groups of even degree and alternating groups of odd degree.

In Theorem \ref{symemtric_odd_degree_and_alternating_even_degree}, the strategy was to look for edges with conjugacy classes of elements having two cycles. This approach is not available anymore. Indeed, in alternating groups of odd degree, elements with two cycles do not exist; and in symmetric groups of even degree, such elements belong to $A_n$, hence one must take care of the parity of elements when dealing with generation. For these elementary reasons, our strategy will be to look for edges with elements having three cycles. This is where we will make use of \cite{GMPS} which, as already mentioned in the introduction, classifies the primitive permutation groups having elements with at most four cycles. 

\begin{theorem}
\label{symmetric_groups_even_degree}
Let $n \geq 8$ be an even integer. Then $\Xi(S_n)$ is connected and $d(\Xi(S_n)) \leq 6$.
\end{theorem}

\begin{proof}
We first assume $n\geq 12$, and consider the remaining cases at the end of the proof. Let $\p$ and $\q$ be two vertices of $\Xi(S_n)$ joined by an edge. One of the two, say $\p$, must correspond to odd permutations. Assume $\p$ is not adjacent to $(i,n-i)$ for any $1 \leqslant i \leqslant n/2$ coprime with $n$. By Lemma \ref{i_coprime_n}, every such $i$, with $i \neq 1$, is partial sum in $\p$. We now show that also $i=1$ is partial sum in $\p$. Assume this is not the case. Since $\p$ is not adjacent to $(1,n-1)$, we deduce from Lemma \ref{intransitive_two_cycles} and Theorem \ref{jones_theorem} that $\p$ is contained in one between $\AGL_d(2), \PGL_2(p), M_{11}, M_{12}, M_{24}$. The last three are excluded because they are subgroups of $A_n$. Assume $\p \in \AGL_d(2)$ or $\PGL_2(p)$. Since every $3 \leqslant i \leqslant n/2$ coprime with $n$ is partial sum, $\p$ must have odd parts: let $a$ be the smallest such part. Since $i=1$ is not partial sum, we have $a \geq 3$. On the other hand, by assumption $\p$ corresponds to odd permutations, hence it has even parts. It follows that $(1^n) \neq \p^a$ fixes a number of points which is greater or equal to $3$, and which is a multiple of $a$. This contradicts Lemma \ref{affine_and_projective}. Therefore, $i=1$ is indeed partial sum in $\p$.

Now we divide the cases $n \equiv 2$ mod $4$ and  $n \equiv 0$ mod $4$.

Assume $n \equiv 2$ mod $4$. Then, $n/2-2$ is coprime with $n$, hence it is partial sum in $\p$. Write $\p=(a_1, \ldots, a_t)$, and assume $n/2-2= \sum_{i=1}^h a_i$. If $a_k=1$ for some $k \in \{1, \ldots, h\}$, then $n/2-3$ is partial sum in $\p$. Otherwise, $n/2-1$ is partial sum in $\p$. We show that $\q$ is adjacent in the first case to $\mathfrak a_1=(1,n/2-3,n/2+2)$, and in the second case to $\mathfrak a_2=(1,n/2-2,n/2+1)$.

By the considerations above, $\q$ and $\mathfrak a_i$ do not share intransitive subgroups. Moreover, $\mathfrak a_1$ and $\mathfrak a_2$ correspond to odd permutations, and belong to no transitive imprimitive subgroups by Lemma \ref{intransitive_three_cycles}. Finally, $\mathfrak a_1$ and $\mathfrak a_2$ belong to no core-free primitive subgroups by \cite[Theorem 1.1]{GMPS}. We have therefore our desired edge between $\q$ and $\mathfrak a_1$ or $\mathfrak a_2$.

Assume now $n \equiv 0$ mod $4$. We employ the same argument as above, with $n/2-2$ replaced by $n/2-1$. The same reasoning  lead us to look for an edge between $\q$ and $\mathfrak b_1=(1,n/2-2,n/2+1)$ or $\mathfrak b_2=(1,n/2,n/2-1)$. Again, $\q$ and $\mathfrak b_i$ do not share intransitive subgroups. It follows from \cite[Theorem 1.1]{GMPS} that $\mathfrak b_1$ and $\mathfrak b_2$ are not contained in core-free primitive subgroups (for $\mathfrak b_2$ we may also use Corollary \ref{cycle_fixing_3_points_power}). Regarding maximal transitive imprimitive subgroups, we only have that $\mathfrak b_2$ is contained in $S_{n/2} \wr S_2$. However, by construction we consider $\mathfrak b_2$ only when $n/2$ is partial sum in $\p$, so that $\p$ belongs to $S_{n/2} \wr S_2$. Since $\p$ and $\q$ are adjacent, we deduce that $\q$ is not contained in $S_{n/2} \wr S_2$. Therefore we have an edge between $\q$ and $\mathfrak b_1$ or $\mathfrak b_2$.

Now we deduce the connectedness of $\Xi(S_n)$ and the bound to the diameter. The considerations above imply that an edge with $\mathfrak a_i$ and $\mathfrak b_i$ concerns only intransitive subgroups (i.e., partial sums), except for $\mathfrak b_2$, where one has to deal also with $S_{n/2} \wr S_2$.

Assume first $n \equiv 2$ mod $4$. The argument given above shows that every vertex of $\Xi(S_n)$ has distance at most $2$ from one among $\mathfrak a_1, \mathfrak a_2$ and $(i,n-i)$ for some $1 \leq i \leq n/2$ coprime with $n$. Hence, in order to conclude it is sufficient to show that these vertices have pairwise distance at most $2$. For $n = 14$, this can be checked directly. Assume then $n > 14$. We show that all these vertices are adjacent to $(2^2,n-4)$, which clearly concludes the proof. For all the vertices except $(1,n-1)$, this follows from Lemma \ref{i_coprime_n} and from the considerations of the previous paragraph. For $(1,n-1)$, by Theorem \ref{jones_theorem} we need to exclude the sharing of $\AGL_d(p), \PGL_2(p), M_{24}$. The last is contained in $A_n$, while $(2^2, n-4)$ is not. Moreover $n \equiv 2$ mod $4$, hence $n$ is not a power of $2$ and we do not have affine subgroups. Finally, $(1^4, n-4)=(2^2,n-4)^2$ fixes at least $4$ points, hence it does not belong to $\PGL_2(p)$ by Lemma \ref{affine_and_projective}(2). This concludes the proof in case $n \equiv 2$ mod $4$.

Assume now $n \equiv 0$ mod $4$. We assume first $n >12$. As in case $n \equiv 2$ mod $4$, in order to conclude it is sufficient to prove that the vertices $\mathfrak b_1, \mathfrak b_2$ and $(i,n-i)$ with $1 \leq i \leq n/2$ coprime with $n$ have pairwise distance at most $2$. The vertices $\mathfrak b_1, \mathfrak b_2$ and $(i,n-i)$ with $1 \leqslant i \leqslant n/2$ coprime with $n$ and $i \neq 3,5$ are adjacent to $(2,3,n-5)$. The vertices $(i,n-i)$ with $1 \leqslant i \leqslant n/2$ coprime with $n$ and $i \neq 1, n/2-1$ are adjacent to both $\mathfrak b_1$ and $\mathfrak b_2$. The vertices $(i,n-i)$ with $2\leq i\leq n/2$ coprime with $n$ are adjacent to $(2^2,n-4)$. By \cite[Theorem 1.1]{GMPS} (which in the affine case relies on \cite[Theorem 1.5]{GMPS1}) we deduce that $(2^2,n-4)$ is not contained in affine subgroups, hence also $(1,n-1)$ is adjacent to $(2^2,n-4)$. These considerations imply indeed that the vertices have pairwise distance at most $2$.

Consider now the case $n=12$. The argument of the previous paragraph does not work, and we need more detailed inspection. Let $\p$ and $\q$ be as at the beginning of the proof, with $\p$ corresponding to odd permutations, and such that $1$ and $5$ are partial sum in $\p$. If $2$ is not partial sum and $4$ is partial sum then it is easy to deduce $\p=(1,4,7)$. If $2$ and $4$ are not partial sums then $\p=(1,5,6)$. Assume now $2$ is partial sum. If $3$ is not partial sum then it is easy to check that $\p$ must have four cycles, false. If $3$ is partial sum and $4$ is not partial sum, then $5$ cannot be partial sum: false. If $4$ is partial sum, then $\p$ is isolated unless $6$ is not partial sum and $\p$ is adjacent to $(12)$. With this more detailed information, it is not difficult to deduce $d(\Xi(S_{12}))\leq 6$. The proof of the theorem for $n\geq 12$ is now concluded. In the next lemma we consider the case $n=8$. The case $n=10$ can be dealt with similarly and we omit the details.
\end{proof}
We compute the exact diameter of $\Xi(S_8)$: this shows that the upper bound in Theorem \ref{main_theorem_2} can be attained.
\begin{lemma}
\label{lemma_graph_S8}
The graph $\Xi(S_8)$ is connected with diameter $6$.
\end{lemma}

\begin{proof}
In Figure \ref{fig:S8} we have drawn the graph $\Xi(S_8)$. The group $S_8$ has $21$ nontrivial conjugacy classes; one can compute explicitly the neighborhood of each of them in $\Lambda(S_8)$. We can save some computations in view of the following observations. Whenever $4$ is a partial sum in a partition $\p$, then $\p$ is not adjacent to partitions having only parts of even length, because of the sharing of $S_4\wr S_2$. It follows that if in a partition $\p$ the integers $1,3$ and $4$ are partial sums, then $\p$ is isolated. This implies that the set of vertices of $\Xi(S_8)$ is a subset of  \[A:=\{(1^3,5),(1^2,6),(1,2,5),(1,7),(3,5),(2,3^2),(4^2),(2,6),(2^4),(2^2,4),(8)\}.
\]

Note also that partitions which, for every odd integer $\ell$, have an even number (possibly zero) of parts of length $\ell$, are not adjacent to partitions having only parts of even length, because of $S_2\wr S_4$. We observe finally that the only core-free maximal primitive subgroups of $S_8$ are $\AGL_3(2)$ and $\PGL_2(7)$ (up to conjugation). Among the partitions in $A$, $\AGL_3(2)$ contains $(1,7), (4^2), (2,6), (2^4)$, and $\PGL_2(7)$ contains $(1^2,6), (1,7), (4^2), (2^4), (8)$. It is now easy to draw the graph.
\end{proof}

\begin{figure}
    \centering
    \begin{tikzpicture}

\node at (0,0) (testa) {$(3,5)$};
\node at (0,1.2) (spalla sinistra) {$(8)$};
\node at (2,1.2) (ginocchio sinistro) {$(1,2,5)$};
\node at (4,1.2) (piede sinistro) {$(4^2)$};
\node at (-2,1.2) (mano sinistra) {$(1^3,5)$};
\node at (2,0) (spalla destra) {$(2^2,4)$};
\node at (4,0) (ginocchio destro) {$(1,7)$};
\node at (6,0) (piede destro) {$(2,3^2)$};
\node at (-2,0) (mano destra) {$(1^2,6)$};

\draw (testa) -- (spalla sinistra);
\draw (testa) -- (mano destra);
\draw (testa) -- (spalla destra);
\draw (spalla destra) -- (ginocchio destro);
\draw (piede destro) -- (ginocchio destro);
\draw (spalla sinistra) -- (ginocchio sinistro);
\draw (spalla sinistra) -- (mano sinistra);
\draw (piede sinistro) -- (ginocchio sinistro);

\end{tikzpicture}
    \caption{The graph $\Xi(S_8)$.}
    \label{fig:S8}
\end{figure}

The last case to consider is alternating groups of odd degree.

\begin{theorem}
\label{alternating_groups_odd_degree}
Let $n \geqslant 5$ be an odd integer. Then, $\Xi(A_n)$ is connected and $d(\Xi(A_n)) \leqslant 6$.
\end{theorem}

\begin{proof}

The cases $n=5,7,9$ can be checked explicitly; we omit the details and assume $n\geq 11$. The cases $n=11,13,17,19$ have been considered in Theorem \ref{lemma_better_estimate_diameter} and in the comments following it. Therefore we need to consider the cases $n=15$ and $n \geq 21$. We first assume $n \geq 21$, and deal with the case $n=15$ at the end of the proof. Throughout the proof, recall Lemmas \ref{intransitive_imprimitive} and \ref{maximal_subgroups_conjugacy}.

Let $x, y \in A_n$, and assume $x^{A_n}$ is adjacent to $y^{A_n}$ in $\Xi(A_n)$. We will show that $x^{A_n}$ or $y^{A_n}$ is adjacent to at least one among $(1^2,n-2), (2^2,n-4)$, every class with cycle type $(1,3,n-4), (1,4,n-5)$ and $(2,8,n-10)$. This will show that every vertex of $\Xi(A_n)$ has distance at most $2$ from one of these vertices. The argument given in the next two paragraphs shows that all these vertices are adjacent to $(6^2,n-12)$. It will follow that every vertex of $\Xi(A_n)$ has distance at most $3$ from $(6^2,n-12)$, which clearly will conclude the proof.

Let us analyze the classes $(1^2,n-2), (2^2,n-4), (1,3,n-4), (1,4,n-5)$ and $(2,8,n-10)$. By Lemma \ref{intransitive_three_cycles}, only $(1,4,n-5)$ and $(2,8,n-10)$ belong to some maximal transitive imprimitive subgroup: they are contained in $S_5 \wr S_{n/5}$. Regarding maximal core-free primitive subgroups, by Theorem \ref{jones_theorem} $(1^2,n-2)$ belongs only to $\PGaL_2(q)$ with $n=q+1$. Moreover, $(2^2,n-4)$, $(2,8,n-10)$ and $(1,3,n-4)$ belong to no core-free primitive subgroups: the first two by Corollary \ref{cycle_fixing_3_points_power}, the last by \cite[Theorem 1.1]{GMPS}. Again by \cite[Theorem 1.1]{GMPS}, $(1,4,n-5)$ belongs possibly only to $\AGL_d(5)$.

Now consider the class $(6^2,n-12)$. The maximal transitive imprimitive subgroups it is contained in are $S_3 \wr S_{n/3}$ and $S_{n/3} \wr S_3$. Moreover, by Lemma \ref{affine_and_projective} it belongs neither to $\AGL_m(p)$ nor to $\PGaL_2(q)$. Therefore, as claimed $(1^2,n-2), (2^2,n-4)$, every class with cycle type $(1,3,n-4), (1,4,n-5)$ and $(2,8,n-10)$ are adjacent to $(6^2,n-12)$.

Hence, in order to conclude the proof it is sufficient to prove the initial claim, i.e., to prove that $x^{A_n}$ or $y^{A_n}$ is adjacent to at least one among $(1^2,n-2), (2^2,n-4)$, every class with cycle type $(1,3,n-4), (1,4,n-5)$ and $(2,8,n-10)$. In a previous paragraph we determined the maximal overgroups of these classes. In the following, we will freely use this information with no further mention.

Denote by $c(x)$ and $c(y)$ the cycle types of $x$ and $y$ respectively. Assume there exists $z \in \{ x,y\}$ such that $1$ and $2$ are not partial sums in $c(z)$; without loss of generality, $z=x$. Then, by Theorem \ref{jones_theorem} either $x^{A_n}$ is adjacent to $(1^2,n-2)$, or $x \in \PGaL_2(q)$. If moreover $4$ is not partial sum in $c(x)$ then $x^{A_n}$ is adjacent to $(2^2,n-4)$. Assume then that $4$ is partial sum in $c(x)$. The unique possibility is $c(x)=(4, \ldots)$, from which $1 \neq x^4$ fixes an even number of points greater than $3$, from which $x \notin \PGaL_2(q)$ by Lemma \ref{affine_and_projective}(2).

Therefore, we assume (without loss of generality) that $1$ is partial sum in $c(x)$ and $2$ is partial sum in $c(y)$, that is, $c(x)=(1, \ldots)$ and $c(y)=(2, \ldots)$. Then, either $x^{A_n}$ is adjacent to $(2^2,n-4)$, or $4$ is partial sum in $c(x)$. In the latter case we have $c(x)=(1,3,\ldots)$ or $c(x)=(1,4,\ldots)$. If $c(x)=(1,3,\ldots)$ then $y^{A_n}$ is adjacent to $(1,3,n-4)$. If $c(x)=(1,4,\ldots)$ then $y^{A_n}$ is adjacent to $(1,4,n-5)$ unless $y$ is contained in $\AGL_d(5)$ or $S_5 \wr S_{n/5}$. The option $y \in \AGL_d(5)$ is excluded by Lemma \ref{affine_and_projective}(1) because $y^2$ fixes $2$ points. If $y \in S_5 \wr S_{n/5}$, since $5$ is partial sum in $c(x)$, the unique possibility for the $2$-cycle is to act as a $2$-cycle on the blocks, from which $c(y)=(2,8,\ldots)$. At this point, since certainly $x \notin S_5 \wr S_{n/5}$, $x^{A_n}$ is adjacent to $(2,8,n-10)$. This concludes the proof in case $n \geq 21$.

There remains the case $n=15$.  Let $x, y \in A_n$, and assume $x^{A_n}$ is adjacent to $y^{A_n}$ in $\Xi(A_n)$. We prove that $x^{A_n}$ or $y^{A_n}$ is adjacent to at least one between $(1^2,13)$, $(1,7^2)$ and $(15)$ (this makes sense by Lemmas \ref{unico_vertice} and \ref{n_cycle}). By Theorem \ref{jones_theorem} and Lemma \ref{intransitive_three_cycles}, $(1^2,13)$ and $(15)$ are adjacent. Moreover $(1,7^2)$ and $(15)$ are both adjacent to $(2^2,11)$; and $(1,7^2)$ and $(1^2,13)$ are both adjacent to $(3^2,9)$. It follows that $(1^2,13)$, $(1,7^2)$ and $(15)$ have pairwise distance at most $2$, from which indeed $d(\Xi(A_{15})\leq 6$.

By Theorem \ref{jones_theorem}, if $x^{A_n}$ and $y^{A_n}$ are not adjacent to $(1^2,13)$ then (without loss of generality) $c(x)=(1,\ldots)$ and $c(y)=(2, \ldots)$. Assume now $x^{A_n}$ and $y^{A_n}$ are not adjacent to $(15)$. We consider the various possibilities. Notice that $x \notin S_3 \wr S_5$, because otherwise $2$ would be partial sum in $c(x)$. Direct inspection (using for instance GAP) shows that if $x \in \PGL_4(2)$ then it must be $c(x)=(1,7^2)$; hence we may assume this is not the case. Therefore, by Theorem \ref{jones_theorem}, it must be $x \in S_5 \wr S_3$, and $y \in S_3 \wr S_5$ or $y \in \PGL_4(2)$. Since $x \in S_5 \wr S_3$, we have $c(x)=(1,4,\ldots)$, hence $4$ is not partial sum in $c(y)$. Inspection immediately implies $y \notin \PGL_4(2)$, from which $y \in S_3 \wr S_5$. In $c(y)$, the $2$-cycle either acts trivially on the blocks, or acts as a $2$-cycle on the blocks. In the first case, $1$ is partial sum in $c(y)$, and in the second case, $4$ is partial sum in $c(y)$. In both cases we get a contradiction, and the proof is finished.
\end{proof}

Now the proof of Theorem \ref{main_theorem_2} and Theorem \ref{main_theorem_improvement_diameter} follows immediately from Theorems \ref{lemma_better_estimate_diameter}, \ref{symemtric_odd_degree_and_alternating_even_degree}, \ref{symmetric_groups_even_degree} and \ref{alternating_groups_odd_degree}. 


\section{Some comments on Conjecture \ref{conjecture_diameter_at_most_4}}
\label{subsection_comments_diameter_4} Conjecture \ref{conjecture_diameter_at_most_4} states that, if $G \in \{A_n,S_n\}$, then up to finitely many exceptions one has $d(\Xi(G)) \leq 4$. Here we reduce this conjecture to the following one (and in fact to something much weaker: see Remark \ref{conjecture_weak}).
\begin{conjecture}
\label{conjecture_number_cycles}
Let $G \in \{A_n,S_n\}$. There exists an absolute constant $c>0$ such that if $x^G$ is a vertex of $\Xi(G)$, then $x^G$ is adjacent to a class which has at most $c$ cycles.
\end{conjecture}

A way to think about this is that, since $x^G$ is a vertex of $\Xi(G)$, by definition $x^G$ is adjacent to some other class $y^G$. It seems conceivable that, summing the parts of the cycle type of $y$ in a suitable way, one obtains that $x^G$ is indeed adjacent to some $z^G$, where $z$ has a bounded number of cycles. In fact, we believe that the value of $c$ should be rather small, say at most $4$: by \cite{GMPS}, only ``few'' core-free primitive subgroups contain elements having at most $4$ cycles. 

We now record a consequence of the Prime Number Theorem.

\begin{theorem}
\label{prime_number_theorem}
Fix $\xi > 0$. Denote by $\pi(n)$ the number of primes less or equal to $n$. Then, $\pi((1+\xi)n)-\pi(n)$ is asymptotic to $\xi n/\ln n$.
\end{theorem}

\begin{proof}
The Prime Number Theorem states that $\pi(n)$ is asymptotic to $n/\ln n$, hence the statement follows from an easy computation.
\end{proof}

\begin{theorem}
Conjecture \ref{conjecture_number_cycles} implies Conjecture \ref{conjecture_diameter_at_most_4}.
\end{theorem}
\begin{proof}
It is sufficient to show that, for $n$ large, vertices which have at most $c$ cycles have pairwise distance at most $2$ in $\Xi(G)$. Let $x^G$ and $y^G$ be two such vertices, and denote by $c(x)$ and $c(y)$ the cycle type of $x$ and $y$, respectively.

We first claim that, if $n$ is sufficiently large, then there exist distinct prime numbers $p$ and $r$ such that:
\begin{itemize}
    \item[(a)] $n/3 < p,r \leq n/2$,
    \item[(b)] $p,r$ and $p+r$ are not partial sum in $c(x)$ and in $c(y)$,
    \item[(c)] $p, r$ and $p+r$ do not divide $n$.
\end{itemize}
Let us prove the claim. By Theorem \ref{prime_number_theorem}, the number of primes contained in the interval $(n/3,n/2]$ is asymptotic $n/6\ln n$. The number of divisors $d(n)$ of $n$ is much smaller: it is known that $d(n)=o(n^{\eps})$ for every fixed $\eps >0$ (cf. \cite[Theorem 13.12, or Exercise 13 p. 303]{Apo}).

Notice now that there are at most $2^{c+1}=O(1)$ integers $i$ such that $i$ is partial sum in at least one between $c(x)$ and $c(y)$.

The claim now follows, just because the number of primes in the interval $(n/3,n/2]$ is much larger than all the other quantities considered above, hence among all possibilities for $p$ and $r$ we certainly find one satisfying (b) and (c).

At this point we conclude the proof. If $G=A_n$ and $n$ is even, or $G=S_n$ and $n$ is odd, then both $x^G$ and $y^G$ are adjacent to $(p,n-p)$. Indeed, by (c) $(p,n-p)$ is not contained in transitive imprimitive subgroups, and a power of $(p,n-p)$ is a $p$-cycle, hence does not lie in core-free primitive subgroups by a classical theorem of Jordan (\cite[Theorem 13.9]{Wie}). Assume now $G=A_n$ with $n$ odd, or $G=S_n$ with $n$ even. Then we claim that both $x^G$ and $y^G$ are adjacent to every class with cycle type $(p,r,n-p-r)$. It is easy to deduce from Lemma \ref{intransitive_three_cycles} and from (c) above that $(p,r,n-p-r)$ does not belong to transitive imprimitive subgroups of $S_n$. Moreover, since $\text{min}\{p,r\} > n/3$ we have $n-p-r < \text{min}\{p,r\}$. It follows that a power of $(p,r,n-p-r)$ is a $p$-cycle, hence does not belong to core-free primitive subgroups by Jordan's theorem.

We have shown that $x^G$ and $y^G$ have distance at most $2$ in $\Xi(G)$, and the proof is concluded.
\end{proof}

\begin{remark}
\label{conjecture_weak}
In the proof of the previous theorem, we have used a much weaker hypothesis than the validity of Conjecture \ref{conjecture_number_cycles}. Indeed, the same argument works provided each vertex $x^G$ of $\Xi(G)$ is adjacent to a vertex $y^G$ such that the number of integers $i$ which are partial sum in (the cycle type of) $y$ is at most $\delta n/\log n$ for some explicit fixed constant $\delta$. This statement seems more suitable for a combinatorial proof. Any proof of this sort would very likely avoid the use of the CFSG.
\end{remark}

We conclude with a lemma providing a lower bound to $d(\Xi(G))$ in some cases.

\begin{lemma}
\label{lower_bound_diameter}
Assume $n \geq 7$ is prime. Then, $d(\Xi(S_n)) \geq 4$.
\end{lemma}

\begin{proof}
The strategy is to define two partitions $\p$ and $\q$ such that $\p$ is adjacent in $\Xi(S_n)$ only to $(n)$, and $\q$ is adjacent in $\Xi(S_n)$ only to $(1,n-1)$. Since $(n)$ and $(1,n-1)$ are not adjacent because of the sharing of $\AGL_1(n)$, this will prove indeed that $d(\Xi(S_n)) \geq 4$.

Define $\p=(1^{(n-1)/2}, (n+1)/2)$ if $n \equiv 3$ mod $4$ and $\p=(1^{(n+1)/2},(n-1)/2)$ if $n \equiv 1$ mod $4$. Note that $\p$ corresponds to odd permutations. Moreover, every $1 \leq i \leq n/2$ is a partial sum in $\p$, hence $\p$ is not adjacent to partitions having at least $2$ parts. Finally, $\p$ is adjacent to $(n)$ by Theorem \ref{jones_theorem}.

Now define $\q=(2^{(n-3)/2},3)$ if $n \equiv 3$ mod $4$ and $\q=(2^{(n-9)/2}, 3^3)$ if $n \equiv 1$ mod $4$. Note that $\q$ corresponds to even permutations, so it is not adjacent to $(n)$. It is easy to check that every $2 \leqslant i \leqslant n/2$ is partial sum in $\q$, hence $\q$ is adjacent to nothing different from $(1,n-1)$. Finally, since $\q \notin \AGL_1(n)$, we deduce by Theorem \ref{jones_theorem} that $\q$ is indeed adjacent to $(1,n-1)$. This concludes the proof of the lemma.
\end{proof}

Similar methods should suffice to give sensible lower bounds to $d(\Xi(G))$ in all the various cases. However, the details would become slightly technical, as one would need to specialize the argument depending on the arithmetic of $n$.

\bibliography{references}
\bibliographystyle{alpha}
\nocite{GAP4}

\end{document}